\documentclass[12pt]{amsart} \usepackage{amscd}

\newcommand{\Aa}{\mathcal{A}}

\newcommand{\End}{\operatorname{End}}
\newcommand{\Ext}{\operatorname{Ext}}

\newcommand{\g}{\mathfrak{g}}

\newcommand{\Hom}{\operatorname{Hom}}
\newcommand{\Id}{\operatorname{Id}}

\newcommand{\ind}{\operatorname{ind}}

\newcommand{\Ker}{\operatorname{Ker}}
\newcommand{\OO}{\mathcal{O}}
\newcommand{\Pic}{\operatorname{Pic}}
\newcommand{\PP}{\mathcal{P}}
\newcommand{\Proj}{\operatorname{Proj}}
\newcommand{\Q}{\mathcal{Q}}
\newcommand{\SL}{\operatorname{SL}}

\newcommand{\slq}{\mathfrak{sl}_q(2)}
\newcommand{\Spec}{\operatorname{Spec}}

\newcommand{\nc}{\mathbb{C}}
\newcommand{\np}{\mathbb{P}}
\newcommand{\nq}{\mathbb{Q}}
\newcommand{\nz}{{\mathbb{Z}}}

\newtheorem{thm}{Theorem}[section]
\newtheorem{lm}[thm]{Lemma}
\newtheorem{prop}[thm]{Proposition}
\newtheorem{crl}[thm]{Corollary}

\theoremstyle{definition}

\newtheorem{df}{Definition}[section]
\newtheorem{ex}{Example}[section]

\theoremstyle{remark}

\newtheorem{rem}{Remark}%[section]
\newtheorem{ack}{Acknowledgment}

\begin{document}

\title
{Semiquantum Geometry}

\author[N. Reshetikhin]
{Nicolai Reshetikhin}
\address
{Department of Mathematics, University of California, Berkeley, CA
94720-3840, USA}
\email{reshetik@math.berkeley.edu}

\author[A. A. Voronov]
{Alexander A. Voronov}
\address
{Department of Mathematics, University of
Pennsylvania, Philadelphia, PA 19104-6395, USA}
\email{voronov@math.upenn.edu}

\thanks{Research of Reshetikhin was supported in part by NSF grant
DMS-9296120, of Voronov by NSF grant DMS-9402076, and of Weinstein by
NSF grant DMS-93-09653 and DOE contract DOE DE-FG03-93ER25177.}

\author[A. Weinstein]
{Alan~Weinstein}
\address
{Department of Mathematics, University of California, Berkeley, CA
94720-3840, USA}
\email{alanw@math.berkeley.edu}

\date{June 5, 1996}

\subjclass{46L87; Secondary 81S10}

\dedicatory{Dedicated to Yu.\ I. Manin on the occasion of his 60th birthday}
\maketitle

In this paper we study associative algebras with a Poisson algebra
structure on the center acting by derivations on the rest of the
algebra. These structures appeared in the study of quantum groups at
roots of 1 and related algebras. They also appeared in the study of
the representation theory of affine Lie algebras at the critical
level. We also believe that these \emph{Poisson fibred algebras} play
an essential role in Fedosov's quantization of matrix-valued functions
on symplectic manifolds. They are also implicit in the work of Emmrich
and one of the authors \cite{em-we:transport} on multicomponent WKB
approximations. In all these cases, this special Poisson algebra
structure is induced by a one parameter family of deformations.

We also take up the general study of
noncommutative spaces which are close to enough commutative ones so
that they contain enough points to have interesting commutative
geometry. One of the most striking uses
of our noncommutative spaces is the quantum Borel-Weil-Bott
Theorem~\ref{BWB} for quantum $\slq$ at a root of unity, which comes
as a calculation of the cohomology of actual sheaves on actual
topological spaces.  The idea of
such noncommutative spaces is not new and has been approached by many
mathematicians, starting with A.~Grothendieck \cite{ega}, who treated
the most general situation of a scheme ringed with a sheaf of
noncommutative algebras.  Later work
\cite{berezin,bl,l,manin:gauge} treated
a particular case, supermanifolds or
$\nz_2$-graded manifolds, which are in some sense quantum manifolds at
$q=-1$, a square root of unity. The idea of unifying supermathematics
with noncommutative geometry was made explicit by Manin, see
\cite{manin:noncomm}.

In a series of papers \cite{a,atv1,av}, M.~Artin, J.~Tate and M.~Van den Bergh
worked out the case of quantum projective planes,
being mostly interested in the ones corresponding to an elliptic curve
in $\np^2$ and an automorphism of the curve. When the automorphism was
of finite order, they were in the situation of the present paper, but
only in dimension two. Yu.~I. Manin \cite{manin:tori} noticed the same
pattern, working with quantum abelian varieties at roots of unity.
B.~Parshall and Jian Pan Wang \cite{pw} made a very comprehensive
study of quantum groups and homogeneous spaces at roots of unity at
the level of Hopf algebras and their comodules. And finally, C.~De
Concini, V.~G.  Kac and C.~Procesi \cite{dkp} studied orbits of the
coadjoint representation of quantum groups at roots of unity and came
across similar spaces.

For a quantum geometric object $X_q$ (a deformation of a classical
object), it would be interesting to consider as a whole the hierarchy
of quantum objects $X_\epsilon$ over roots $\epsilon$ of unity
of all possible orders. The analogy with Diophantine
geometry beyond the Frobenius morphism (which in this work relates
$X_\epsilon$ with $X_1$), that is, for instance, studying the
corresponding Galois groups and relating them to
$\operatorname{Gal}(\bar \nq/\nq)$, defining the corresponding
$\zeta$-functions, etc., seems to be very promising, both for
noncommutative geometry and number theory.

\begin{ack}
We are grateful to M.~Gerstenhaber, A.~Givental, V.~Kac, Yu.~I. Manin,
M.~Movshev, A.~S. Schwarz, J.~Stasheff, and A. Vaintrob for helpful suggestions
and interesting discussions.  The first author thanks the Miller Institute for
partial support, and the second author thanks the Max Planck Institute for
Mathematics at Bonn for hospitality during the summer of 1995, when large part
of the work was done.
\end{ack}

\section{Poisson fibred algebras}
\label{quantum}

\subsection{} Here and below we assume that all objects are defined over
a field $\kappa$ of characteristic zero.

As usual, by a Poisson algebra we mean a commutative algebra with a
Lie algebra structure on it (given by a Poisson bracket), such that
the Poisson bracket acts on the commutative algebra by derivations.  A
Poisson manifold is a smooth manifold $X$ with a Poisson algebra
structure on the smooth functions on $X$ under pointwise
multiplication.

\begin{df} A \emph{Poisson fibred algebra} is a triple $(A,Z,\{\cdot,\cdot\})$,
where $A$ is an
associative algebra with a unit 1, $Z$ is its center,
and $\{\cdot,\cdot\}$ is a bracket
\begin{align*}
Z \times A & \to A,\\
(f,a) & \mapsto \{f,a\}
\end{align*}
which has the following properties.
\begin{enumerate}
\item Its restriction to $Z \times Z$ provides $Z$ with the structure of a
Poisson algebra.
\item It gives an action by derivations of this Poisson algebra on
the whole algebra $A$, such that
\begin{align*}
\{f, a b\}  &  = \{f,a\} b + a \{f,b\},\\
\{fg, a\} & = f\{ g, a \} + g \{ f, a\},
\end{align*}
for $f, g \in Z$, $a, b \in A$.
\end{enumerate}
\end{df}

Notice that the last two equations imply
\[
\{f,1\} =  \{1,a\} = 0 \qquad \text{for all $f \in Z$, $a \in A$}.
\]

\begin{df} The obstruction $\Phi(f,g) \in Der(A)$
to $A$ being a Lie algebra module over $Z$
\begin{equation} \label{curvat}
\Phi(f,g) (a) = \{ \{ f, g \}, a \}- \{ f, \{ g, a \}\} + \{ g, \{ f, a\}\},
\end{equation}
where $f,g \in Z$, $a \in A$, is called the \emph{curvature of the
Poisson fibred algebra}.
\end{df}

It is obviously $Z$-linear in $f$, $g$ and $a$ and
skew in $f$ and $g$.

Notice that any linear map $\theta: A\to A$ determines a new bracket
\begin{equation}  \label{gauge}
\{f,g\}_\theta=\{ f, g \}+[\theta(f),g],
\end{equation}
which also provides a structure of Poisson fibred
algebra on $A$.

We will call two Poisson fibred algebras \emph{equivalent}
if their brackets are
related by a transformation
\eqref{gauge}.  The curvature $\Phi(f,g)_\theta$ of the Poisson fibred
algebra with the bracket \eqref{gauge} is related to \eqref{curvat} as
follows
\begin{multline} \label{gaugecur}
\Phi(f,g)_\theta (a)=\Phi(f,g) (a) \\
+[\theta(\{f,g\}),a]
-[\{f,\theta(g)\},a]+[\{g,\theta(f)\},a]-[[\theta(f),\theta(g)],a].
\end{multline}

Thus, strictly speaking, by the curvature of a
 Poisson fibred algebra we should understand the
equivalence class of $\Phi$ under the transformations
\eqref{gaugecur}. In particular, one can say that a
Poisson fibred algebra is flat if its curvature
can be removed by a ``gauge'' transformation \eqref{gauge}.

\subsection{The geometric interpretation of Poisson fibred algebras}

Poisson fibred algebras arise from the following geometric situation.

Recall that a \emph{Lie algebroid} is a vector bundle $L$
over a smooth manifold $X$ together with the structure $[,]$ of
a Lie algebra on its smooth sections and a vector-bundle
morphism $\rho: L  \to TX$ (called the \emph{anchor} of the Lie
algebroid) satisfying the following conditions:
\begin{gather}
[\rho(a), \rho (b)]  = \rho ([a, b]),\notag\\
\label{Leibniz}
[a, f b] = f [a, b] + (\rho(a)f) b,
\end{gather}
where $a$ and $b$ are sections of $L$ and $f$ is a
function on $X$.  For more details
about Lie algebroids, see \cite{kenz}.

Suppose that we have a Poisson manifold $X$. The Poisson structure on
$X$ defines the structure of a Lie algebroid on its cotangent
bundle.  The anchor is defined as contraction $\iota$ with
the bivector field $P \in \Gamma (X, \Lambda^2 T X)$ defining the
Poisson structure:
\begin{align*}
\rho: T^{*} X & \to TX,\\
\rho(\omega) & = \iota(\omega) P = \omega (P).
\end{align*}
The bracket on 1-forms is defined as the one extending the bracket
\[
[df, dg] = d\{f,g\}
\]
by the Leibniz rule \eqref{Leibniz}.  This yields the following
formula for arbitrary 1-forms:
\[
[\omega, \varphi] = L_{\rho (\omega)} \varphi - L_{\rho(\varphi)}
\omega + d \iota(\omega \wedge \varphi) P.
\]

Thinking of a Lie algebroid as a generalization of the tangent bundle,
we would like to define the notion of a Lie algebroid connection in a
vector bundle $V$ over $X$. Consider the \emph{Atiyah algebra}
$\mathcal A (V)$ of $V$, which is the bundle of first-order
differential operators on $V$ with symbol $\sigma = \Id_V \otimes
\xi$, where $\xi$ is a vector field on $X$. There is a natural short
exact sequence
\begin{equation*}
0 \to \End (V) \to \mathcal A (V) \overset{\sigma}{\to} TX \to 0
\end{equation*}
of vector bundles. The Atiyah algebra is a Lie algebroid with respect to the
commutator of differential operators and the anchor map $\sigma$.

If $L$ is a Lie algebroid over $X$, an
anchor preserving vector bundle mapping
\begin{equation}
\label{nabla}
\nabla: L \to \mathcal A (V)
\end{equation}
will be called an $L$-{\em connection on} $V$; we define its
\emph{curvature} as the section $\Phi_\nabla$ of the bundle $\Lambda^2
L^{*} \otimes \Aa (V)$ by the formula
\[
\Phi_\nabla (\omega \wedge \varphi) = \nabla [\omega, \varphi] -
[ \nabla \omega, \nabla \varphi].
\]
An $L$-connection $\nabla$ is a morphism of
Lie algebroids if and only if its curvature vanishes. In this case,
the ``flat connection'' $\nabla$ is also called a \emph{representation of the
Lie algebroid
$L$ in $V$}.  (When $X$ is a point, $L$ is just a Lie algebra, and
this is the usual definition of a representation.)

Now assume
that $V$ is a bundle of associative algebras over a Poisson manifold
$X$ and that $\nabla $ is a $T^{*}X$-connection on $V$ such that
\begin{equation}
\label{der}
\nabla_\omega (ab) = (\nabla_\omega a) b + a (\nabla_\omega b)
\end{equation}
for all 1-forms $\omega $ and $a, b \in C^\infty (X, V)$.

\begin{thm}
\label{main}
The bracket
\[
\{f, a\} = \nabla_{df} a
\]
where $f$ is a function on $X$ and $a$ is a section of $V$, defines
the structure of a Poisson fibred algebra on sections of
$V$. Conversely, let $Z = C^\infty (X)$ be the algebra of smooth
functions on a manifold $X$ and $A = C^\infty (X,V)$ the $Z$-module of
smooth sections of a vector bundle $V$ over $X$. Then any local
structure of a Poisson fibred algebra on the pair $(A,Z)$ gives rise
to a $T^{*}X$-connection $\nabla$ on $V$ satisfying \eqref{der}.  The curvature
of this
Poisson fibred algebra is related to the curvature of $\nabla$ by the
formula
\[
\Phi(f,g) = \Phi_\nabla (df,dg).
\]
\end{thm}

\begin{proof}
The fact that a $T^{*}X$-connection satisfying
\eqref{der} gives a Poisson fibred algebra is a direct verification
of axioms. The converse statement is proved by defining a $T^{*}X$-connection
by the formula
\[
\nabla_\omega a (x) = \{f,a\} (x),
\]
where $\omega$ is a 1-form, $a$ is a section of $V$ and $f$ is a
function on $X$, such that $df(x) = \omega(x)$ at a point $x \in
X$. The fact that for a fixed section $a$ of $V$, the bracket
$\{f,a\}$ is a first-order differential operator on $X$ ensures that
$\{f,a\}(x)$ depends only on the first jet of $f$ at $x \in X$, which
guarantees correctness of our definition of $\nabla$.
\end{proof}

\subsection{Formal deformations}
\label{deform}

Poisson fibred algebras appear naturally in the context of formal
deformations.

Consider an associative algebra $A$ with the center $Z$.
Let $A_h$ be a formal deformation of $A$
such that $A_h$ is a $\kappa[[h]]$-algebra isomorphic to $A[[h]]$ as a
module over $\kappa[[h]]$. Choose such an isomorphism $\phi: A[[h]]
\to A_h$. Let $\phi(1) = 1$. In the sequel we will identify $A[[h]]$
with $A_h$ via this isomorphism.

After this identification the multiplication in $A_h$, sometimes called
$*$-multiplication, is given by a
formal power series
\begin{equation}
f * g = fg + h B_1(f,g) + h^2 B_2(f,g) + \dots.
\label{pair}
\end{equation}
Here $B_i (f,g)$ are $A$-valued bilinear forms on $A$.
The following proposition follows from identities for commutators in
the associative algebra $A_{h}$.

\begin{prop}
\begin{enumerate}
\item The bracket
\begin{equation} \label{br}
\{a, b\}=B_1(a,b) - B_1(b,a)
\end{equation}
determines a Poisson fibred  algebra structure on $(A,Z)$.
\item The curvature of this Poisson fibred algebra is
given by
\[
\Phi(f,g) (a) = [a, B_2(f,g) - B_2 (g,f)].
\]
\end{enumerate}
\end{prop}

If we compose the identification map $\phi: A[[h]] \to A_h$ with any
$\kappa[[h]]$-linear map $\psi: A[[h]] \to A[[h]]$ and assume that
$\psi= \Id +h \psi_1+ \dots$, the bracket \eqref{br} will change to
\begin{equation*}
\{a,b\}' = \{a,b\} + [\psi_1(a),b].
\end{equation*}

The curvature will change according to \eqref{gaugecur} with $\theta
= \psi_1$.

\section{Some examples}

\begin{ex} Any Poisson algebra is a flat Poisson fibred algebra with $A=Z$.

\end{ex}

%\begin{sloppypar}
\begin{ex}
Consider a two dimensional quantum torus.
For any $q\in {\mathbb C}$, the two dimensional
quantum torus is the associative algebra $A_q$ generated by two
invertible elements $u$ and $v$ subject to the relations
\begin{equation}
uv=qvu
\end{equation}
When $q=\epsilon$, where $\epsilon^l=1$, the elements
$u^l$ and $v^l$ and their inverses generate
the center of the algebra $A_\epsilon$, and the algebra
$A_\epsilon$ becomes a Poisson fibred algebra.

The Poisson bracket between elements  $u^l$ and $v^l$ can be easily
computed from \eqref{br}. If we choose the isomorphism
between vector spaces $A_\epsilon$ and $A_q$ identifying
monomials $u^i v^j$ in both algebras, we obtain the following brackets:
\begin{equation}
\{u^l, v^l \}=l^2u^lv^l
\end{equation}
It is easy to check that in this case $B_2(f,g)$ is central if $f$
and $g$ are central. Therefore the Poisson fibred
algebra $A_\epsilon$ has zero
curvature.
\end{ex}

\begin{ex}[\cite{dk}]
\label{kac}
Let $\kappa = \nc$ and $q$ be a formal parameter.  The quantum universal
enveloping algebra of ${\mathfrak sl}_2$ is the associative algebra
$U_q({\mathfrak sl}_2)$ over $\nc$ with generators $E, F, K$ and $K^{-1}$
subject to the following relations:
\begin{align*}
KEK^{-1} & = q^2 E,\\
KFK^{-1} & = q^{-2} F,\\
[E,F] & = \frac{K-K^{-1}}{q-q^{-1}} .
\end{align*}
It has a Hopf algebra structure reflecting its group properties. Here
we will focus on its associative algebra properties when $q$ is a root
of unity.

It is well known that the center of $U_q({\mathfrak sl}_2)$ for a
generic $q$ is generated by the element
\[
c=EF+ \frac{Kq^{-1}+K^{-1}q}{(q-q^{-1})^2} .
\]

Let $l> 2$ be an odd integer.  Specify $q = \epsilon$, a primitive
$l$th root of unity.  According to Section~\ref{deform}, the algebra
$A = U_\epsilon = U_{\epsilon}({\mathfrak sl}_2)$ is an example of a
Poisson fibred algebra. Here we will describe a slight modifictaion
(with $Z$ being a certain central subalgebra of $A$, as opposed to the
whole center) of this algebra explicitly.  The case of an even root
of $1$ is technically a little more complicated, but absolutely
similar.

The center of $U_{\epsilon}({\mathfrak sl}_2)$ is generated by the
elements $E^l, F^l, K^l, K^{-l},c$.  Consider the subalgebra $\nc
[E^l, F^l, K^l, K^{-l}] \subset U_q$. At $q=\epsilon$ it degenerates
to a subalgebra $Z$, which lies in the center of $A = U_\epsilon$, and
$U_\epsilon$ becomes a free $Z$-module of rank $l^3$. The algebra $Z$
has a canonical Poisson bracket
\[
\{a,b\} = \lim_{q \to \epsilon} \frac{\tilde a \tilde b - \tilde b \tilde
a} {q-\epsilon},
\]
where $a, b \in Z$ and $\tilde a , \tilde b \in U_q$ are liftings of
these elements: $a = \tilde a \pmod{q-\epsilon}$, $b = \tilde b
\pmod{q-\epsilon}$.  It also acts by derivations on $U_\epsilon$:
\[
\{a,u\} = \lim_{q \to \epsilon} \frac{\tilde a \tilde u - \tilde u \tilde
a} {q-\epsilon},
\]
but the bracket $\{a,u\}$ is defined up to the addition of $vu -uv$ for
some $v \in U_\epsilon$, because of the ambiguity of lifting of $a$
modulo $q-\epsilon$.

The monomials $F^i K^j E^k$, $i,j,k \in \nz$, $i,k \ge 0$, form a
basis of $U_q$, and this choice of basis identifies $U_\epsilon$ with
a subspace of $U_q$, thereby defining the liftings $\tilde a$,
etc. This allows to compute explicitly the Poisson structure on $X =
\nc^2 \times \nc^*$ in coordinates $x=E^l$, $y=F^l$ and $z=K^l$ and
the action of $Z$ on the vector bundle $V$ over $X$ associated with
the free $Z$-module $A$ in coordinates $E$, $F$ and $K$ on $A$. Here
we use explicit formulas from De Concini-Kac \cite{dk} for commutators
of elements in $A$.
\begin{align*}
\{x,y\} & = d c^{-l} (z
-z^{-1}),\\
\{x, z\} & = - d c^l x
z,\\
\{y, z\} & = d
c^l yz,\\
\{x,E\} & = 0,\\
\{x, F\} & =  (d
c^{l-2}/l) (K \epsilon - K^{-1} \epsilon^{-1}) E^{l-1},\\
\{x, K\} & = - (d c^{l}/l)  x K,\\
\{z, E\} & = (d c^{l}/l) zE,\\
\{z,K\} & = 0,
\end{align*}
etc., where
\begin{gather*}
d = \lim_{q \to \epsilon} \frac{[l]!}{q-\epsilon}\\
[l]! = [l] [l-1] \dots [1] = \frac{q^l - q^{-l}}{q-q^{-1}} \frac{q^{l-1} -
q^{1-l}}{q-q^{-1}} \dots \frac{q - q^{-1}}{q-q^{-1}},\\
c = \epsilon -\epsilon^{-1}.
\end{gather*}
Thus the Poisson tensor on $X$ is equal to
\begin{multline*}
P =
d c^{-l} (z
-z^{-1}) \; \partial/\partial x \wedge \partial/\partial y \\
- d c^l x
z \; \partial/\partial x \wedge \partial/\partial z + d
c^l yz \; \partial/\partial y \wedge \partial/\partial z.
\end{multline*}

The curvature is also computable:
\begin{equation*}
\Phi(x,y) \, a = [\,a, \;\sum_{j=1}^{l-1} \frac{d^2}{([l-j]!)^2 [j]!}
F^{l-j} \prod_{r=j-2l+1}^{2j-2l} [K;r] E^{l-j} \,],
\end{equation*}
where $q$ is specialized to $\epsilon$ in $[l-j]!$ and $[j]!$ and
\begin{equation*}
[K;r] = \frac{K\epsilon^r - K^{-1} \epsilon^{-r}}{\epsilon - \epsilon^{-1}}.
\end{equation*}
In particular, it is clear that the curvature is nonzero. However
\[
\Phi(x,z) = \Phi (y,z) = 0.
\]
\end{ex}
%\end{sloppypar}

\begin{ex}[Quantum tori and abelian varieties, Manin \cite{manin:tori}]
\begin{sloppypar}
A \emph{quan\-tum torus} $T(H,\alpha)$ is a pair consisting of a free
finitely generated abelian group $H$ and a multiplicatively
skew-symmetric bi\-linear form $\alpha: H \times H \to \nc^*$, i.e.,
$\alpha (\chi, \eta) = \alpha (\eta, \chi)^{-1}$, $\alpha (\chi+\chi',
\eta)= \alpha (\chi, \eta) \cdot \alpha (\chi', \eta)$.  The {\it
quantum function ring of the torus} is the vector space $A(H,\alpha)$
freely generated over $\nc$ by symbols $e_\chi$, $\chi \in H$, with
the associative multiplication law
\[
e_\chi e_\eta = \alpha(\chi, \eta) e_{\chi + \eta}.
\]
which equips $A(H,\alpha)$ with the structure of an associative
algebra.
\end{sloppypar}

When $\alpha \equiv 1$, we get the commutative multiplication law for
functions on the usual complex torus $T(H, 1)$, i.e., $\Hom(H,
\nc^*)$. Otherwise it is deformed, with $\alpha$ being thereby a
quantization (multi)parameter.

Let us now introduce the structure of a Poisson fibred algebra on a quantum
torus.

Suppose we came to the root of unity case by deforming parameters of the
quantum torus, that is, deforming the form $\alpha$ along a tangent vector
$\gamma$ to the parameter space, at a point where $\alpha $ takes values in
roots of unity. The tangent space to the parameter space at $\alpha$ can be
identified with the space of skew-symmetric bilinear forms (in the usual
additive sense) $\gamma$ on $H$: a tangent vector $\gamma$ represents the 1-jet
$\alpha (1 + t \gamma + \dots )$ of a curve.

Then we can define a Poisson fibred algebra structure at a root of
unity as follows:
\begin{multline*}
\{ e_\chi, e_\eta \} \\
\begin{split}
& = D_\gamma [e_\chi, e_\eta ]^\sim \\
& =
\lim_{t \to 0} (\alpha(\chi, \eta)(1 + t \gamma (\chi, \eta)) - \alpha^{-1}
(\chi, \eta) (1 + t
\gamma(\chi,\eta))^{-1}) e_{\chi + \eta} /t \\
& = 2 \gamma (\chi, \eta) e_{\chi + \eta},
\end{split}
\end{multline*}
where $\chi, \eta \in H$, at least one of them being in $H' = \Ker \alpha$,
$D_\gamma$ is the directional derivative, and $[,]^\sim$ is the commutator
evaluated at a point $\alpha(1 + t \gamma + \dots)$. When $\chi$ and $\eta$ are
in $H'$, we get the structure of a Poisson algebra on $A(H',1)$, implying a
Poisson structure on the torus $T(H',1)$, and when $\chi \in H'$ and $\eta \in
H$, we get an action of $A(H',1)$ by derivations on $A(H, \alpha)$, which
defines the structure of a Poisson fibred algebra on $(Z,A) = (A(H',1),
A(H,\alpha))$, according to the generalities of Section~\ref{deform}. The
curvature of this Poisson fibred algebra vanishes, because away from the root
of unity point, the pair $(A(H',1), A(H,\alpha))$ deforms as a pair, i.e.,
$A(H',1)$ lifts to a \emph{subalgebra}, $A(H',\alpha |_{H'})$ of $A(H, \alpha)$
for generic $\alpha$. This implies $B_2(f,g) \in A(H',1)$ and therefore
$\Phi(f,g) = 0$, see Section~\ref{deform}.
\end{ex}

\begin{rem}
It would be interesting to study the similar structure on a quantum abelian
variety, where there are not enough global functions and their replacement, the
quantum theta functions do not have a natural associative algebra structure. To
make the problem more concrete, we recall Manin's construction
\cite{manin:tori} of quantum abelian varieties.

The complex torus $T(H,1)$ acts on the quantum torus $T(H,
\alpha)$. More precisely, $T(H,1)$ acts by ring homomorphisms of the
function ring: $f^*: A(H, \alpha) \to A(H, \alpha)$, $f \in A(H, 1) =
\Hom (H, \nc^*)$, by
\[
(f^* e_\chi) = f(\chi) e_\chi.
\]
If $B \subset T(H,1)$ is a subgroup (a period subgroup), then a {\it
quantized theta function} is a formal series $\theta = \sum_\chi
a_\chi e_\chi$ which is automorphic with respect to $B$ with
``linear'' multiplicators $\lambda e_{\chi}$:
\[
b^* \theta = \lambda_b e_{\chi_b} \theta ,
\]
$\lambda_b \ne 0$. In the classical $\alpha \equiv 1$ theory, when a
type $L$ of multiplicators is fixed, theta functions are formal sections
$\Gamma (L)$ of the corresponding line bundle $L$ on the abelian
variety $T(H,1) /B$. When $L$ is a polarization, i.e., a positive line
bundle, the graded ring $\bigoplus_{n\ge 0} \Gamma (L^n)$ is the ring
of homogeneous functions on the abelian variety. Therefore, we can
think of the collection of quantum theta functions as defining a {\it
quantum abelian variety}.

Now suppose the form $\alpha$ takes values at roots of unity. Then we get a
large central subalgebra $Z = A(H',1)$, where $H' = \Ker \alpha$ (the center is
$A(\Ker \alpha^2)$), of the algebra $A(H,\alpha)$. Thus, the latter forms a
bundle $V$ of noncommutative algebras over the usual torus $T(H',1) = \Spec Z$.
Any type of multiplicators such that $\chi_b \in H'$ for $b \in B$ determines a
line bundle $L_0$ on the abelian variety $\mathcal{A} = T(H',1)/B$. If in
addition $L_0$ is a polarization and $L$ an arbitrary type of quantum
multiplicators, then the quantum theta functions of the types $L \otimes
L_0^n$, $n \ge 0$, form a graded module over the homogeneous ring
$\bigoplus_{n\ge 0} \Gamma (\mathcal{A}, L_0^n)$ and therefore define a
coherent sheaf of quantum theta functions on the classical abelian variety
$\mathcal{A}$.

It would be interesting to study related Poisson structures on the
abelian variety $\mathcal{A}$ and the sheaf of quantum theta functions
on it.
\end{rem}

\section{Poisson modules}

The goal of this section is to introduce notions of modules
over Poisson algebras.
The idea is to take the
``infinitesimal part'' of the corresponding notions for
formal deformations.

A Poisson fibred algebra with zero curvature carries the important
algebraic structure of a Poisson module. While Poisson fibred algebras
arise under deformation of an algebra $A$ and a central subalgebra $Z$
of it, Poisson modules show up by deformation of a pair $(A,Z)$, where
$Z$ is a commutative algebra and $A$ a module over it. Geometrically,
this corresponds to the ``semiclassical limit'' of a quantization of a
manifold along with a vector bundle over it.

\subsection{Left and right modules over formal deformations}

Let $A_{h}$ be a formal deformation of a Poisson algebra $A$.
We assume that this is a torsion
free deformation and will identify $A_h$ with $A[[h]]$ as vector
spaces.  Let $M_h$ be a module over $A_h$ which is a formal
deformation of the module $M$ over the commutative algebra $A$
($M=M_h/hM_h$).  Again we assume that $M_h$ is isomorphic to $M[[h]]$
as a vector space and fix this isomorphism.  Then the multiplication
in $A_h$ and the module structure are given by formal power series:
\begin{align*}
f * g  & =  f g + (1/2)\{f, g \}h + B(f,g) h^2 + \dots, \\
f * m  & =  f m + (1/2)\{f, m \} h + \tilde B(f, m) h^2 + \dots
\end{align*}
\begin{rem}
In physical applications, we would work over $\nc$ and let
$h$ be $i$ times Planck's constant.
\end{rem}

Assume that the correspondence $\sigma: h \mapsto -h$
defines an anti-involu\-tion
\begin{equation} \label{invol}
\sigma(f*g) = \sigma(g) * \sigma (f)
\end{equation}
of the associative algebra $A[[h]]$. This means in particular that the
bracket $\{f,g\}$ is skew and $B(f,g)$ is symmetric.

It is obvious that for $f \in A[[h]]$ and $m \in M[[h]]$ the formula
\begin{equation}
\label{right}
m * f := \sigma(\sigma(f) * \sigma(m))
\end{equation}
defines on $M[[h]]]$ the structure of a right $A[[h]]$-module.
For $f \in A$ and $m \in M$, we have
\begin{equation}
\label{right1}
m * f = f m - (1/2)\{f, m \} h + \tilde B(f, m) h^2 - \dots
\end{equation}

This property of formal deformations generalizes the fact that any
left module over a commutative algebra is also a right module over
this algebra, and that any left module over Lie algebra carries a
right module structure as well.

If we follow the ``deformation philosophy'', it is natural to ask the
question: which modules over a commutative algebra $A$ can be deformed
to modules over formal deformations of $A$?  This question suggests
the following definition.

\begin{df}
\label{left}
A vector space $M$ is called a \emph{left Poisson module} over a
Poisson algebra $A$ if a scalar multiplication $fm$ and a bracket $\{f,m\}$
between any elements $f \in A$ and $m \in M$ are defined, such that
\begin{enumerate}
\item  the scalar multipication defines on $M$ the structure of a module over
the commutative
associative algebra $A$;
\item $\{f,g\} m  = f \{g,m\} - \{fg, m\} + \{f, gm\}$;
\item there exists $\tilde B \in \Hom (A \otimes M, M)$, such that
\begin{multline*}
2\{\{f,g\},m\} - \{f,\{g,m\}\} + \{g,\{f,m\}\} \\ = 4(f \tilde B (g,m) -
\tilde B (g,fm) + \tilde B (f,gm) - g \tilde B (f,m)) .
\end{multline*}
\end{enumerate}
\end{df}

\begin{rem}
If we regard $\Hom(M,M)$ naturally as a module over the associative
algebra $A$ and consider the corresponding Hochschild complex
$\Hom(A^{\otimes \bullet}, \Hom (M,M))$ with a differential $d$, then
Property 2 of the above definition is equivalent to saying that the
two-cocycle $\{f,g\}m$ is exact:
\begin{equation} \label{ident-1}
\{f,g\}m = (d \alpha) (f\otimes g,m), \quad \text{with $\alpha (f,m) =
\{f,m\}$}.
\end{equation}
Likewise, Property 3 is equivalent to the existence of a one-cocycle $\tilde
B$, such that
\begin{equation} \label{ident-2}
2\{\{f, g\},m\} - \{f,\{g,m\}\} + \{g,\{f,m\}\} = 4 (d \tilde B ) (f
\otimes g - g \otimes f, m).
\end{equation}
\end{rem}

Unfortunately, we do not have any good geometric interpretation of
Properties 2 and 3 of left Poisson modules.

  From the deformational point of view, left Poisson modules admit first
order deformations extendable to the second order as left associative
modules over formal deformations of $A$.

Similarly \emph{right Poisson modules} are those which admit second
order deformations as right modules over associative formal
deformation of $A$.

\subsection{Poisson modules}

A vector space $M_h$ over ${\mathbb C}[[h]]$ is a bimodule over $A_h$
if it has a structure of left and right modules over $A_h$ and these
two actions commute.

Since any right module over a formal deformation $A_h$ of a Poisson
algebra (provided that $A_h$ has the property \eqref{invol}) is also a
left module, the bimodule structure on $M_h$ is equivalent to two
commuting left module structures.

The special class of bimodules emerges when we assume that the left
and right actions are equal modulo $h$:
\[
f*m-m*f =h \{f,m\}'+O(h^2) \ .
\]
Let $(1/2)\{f,m\}$ and $(1/2)\{m,f\}$ be the first order terms of the
left and right actions respectively. Then
\[
\{f,m\}' = (1/2)(\{f,m\}-\{m,f\}) \ .
\]
The brackets $\{f,m\}'$ satisfy the following identities:
\begin{align}
\{\{f,g\},m\}' &=& \{f,\{g,m\}'\}' -
 \{g,\{f,m\}'\}', \notag \\
\label{ident-3}
\{f,gm\}' & = & g\{f,m\}'+\{f,g\} m, \\
\{fg,m\}' & = & f\{g,m\}'+g\{f,m\}'. \notag
\end{align}
This is why the following definition seems natural.

\begin{df} A vector space $M$ is called a {\it Poisson module}
over Poisson algebra $A$ if $M$ is a module over a commutative algebra
$A$ and if a bilinear
bracket $\{\, , \, \}': A\otimes M\to M$ which satisfies \eqref{ident-3} is
given.
\end{df}

\begin{prop}
Any formal deformation of a pair $(A,M)$ to a left $A[[h]]$-module
$M[[h]]$ with $\sigma : h \mapsto -h$ defining an anti-involution on
$A[[h]]$ and such that the right $A[[h]]$-module structure on $M[[h]]$
defined by \eqref{right} makes it a bimodule defines the structure of
a Poisson algebra on $A$ and a Poisson module on $M$ over $A$ with respect to
the bracket $\{ \, , \, \}$ of \eqref{right1}.
\end{prop}
Under the assumptions of this proposition, $\{ \, , \, \} = \{ \, , \,
\}'$, and therefore, it follows from the properties of $\{ \, , \,
\}'$.

Here is one more argument why this definition is ``natural''.
Poisson modules describe Poisson extensions of Poisson algebras (as
bimodules over associative algebras can be used to construct
associative extensions).  The proof of the next proposition is left to
the reader.

\begin{prop}
Given a Poisson algebra $A$ and a vector space $M$, the structure
of a Poisson module on $M$ is equivalent to the structure of a Poisson
algebra on $A \oplus M$ extending that on $A$ and such that $M \cdot M = \{ M,
M\} = 0$.
\end{prop}

Poisson modules also have a nice operadic meaning, cf.\
\cite{gk,loday}.  Consider the Poisson operad, that is the operad
whose $n$th component $\PP(n)$ is the vector space generated by
elements of degree one in each variable $X_1, \dots , X_n$ of the free
Poisson algebra generated by $X_1, \dots , X_n$. $\{\PP(n)\; | \; n
\ge 1 \}$ is a quadratic operad with two generators $X_1 \cdot X_2$
and $\{ X_1, X_2 \}$, such that $X_1 \cdot X_2 = X_2 \cdot X_1$ and $
\{ X_1, X_2 \} = - \{ X_2, X_1 \}$ and relations
\begin{gather*}
X_1 (X_2 X_3) = ( X_1 X_2 ) X_3,\\
(X_1 X_3) X_2 = X_1 (X_3 X_2),\\
\{ \{ X_1, X_2 \}, X_3 \} + \{ \{ X_2, X_3 \}, X_1 \}
+ \{ \{ X_3, X_1 \}, X_2 \} = 0,\\
\{ X_1, X_2 X_3 \} = \{ X_1, X_2 \} X_3 + X_2 \{ X_1, X_3 \},\\
\{ X_2, X_1 X_3 \} = \{ X_2, X_1 \} X_3 + X_1 \{ X_2, X_3 \},\\
\{ X_3, X_1 X_2 \} = \{ X_3, X_1 \} X_2 + X_1 \{ X_3, X_2 \} .
\end{gather*}
A Poisson algebra is nothing but an algebra over the operad $\PP(n)$
and a Poisson module over a Poisson algebra is exactly a module over
an algebra over the Poisson operad.

\section{Noncommutative spaces of finite type}

Considerations of Section~\ref{quantum} arise when a usual quantum
space, e.g., a quantum group $U_q (\g)$, gets a large center at a
certain value of the parameter $q$. This is how the Poisson structures
enter the picture. On the other hand, even the undeformed
part of that picture, an associative algebra of finite
type over a central subalgebra, is interesting. By analogy with
supermanifolds, we would like to consider a more general {\it
noncommutative space of finite type}, which is either a smooth
manifold with a finite rank bundle of associative algebras or a pair
$(X, \Q_X)$, where $X$ is a scheme with a structure sheaf $\OO_X$ and
$\Q_X$ is a sheaf of $\OO_X$-algebras coherent as a sheaf of
$\OO_X$-modules. Depending on the context (smooth or algebraic), we
will use either kind of noncommutative spaces. From this section on,
we will leave the ``semiclassical'' Poisson world for the
``semiquantum'' world of noncommutative spaces of finite type.

In the algebraic context, every algebra which is finite-dimensional over a
central subalgebra is an example:
\begin{lm}
\label{lem}
Let $A$ be an algebra, $Z$ a subalgebra of its center and $A$ be
finitely generated as a $Z$-module. Then there exists a coherent sheaf
of noncommutative algebras over $\Spec Z$, such that the stalk of the
sheaf over a point $\mathfrak p \in \Spec Z$ is exactly the
localization $A_{\mathfrak p}$ at the corresponding prime ideal
$\mathfrak p$ of $Z$.
\end{lm}
\begin{proof}
The proof just repeats the commutative construction of the structure
sheaf on $\Spec Z$. For every element $f$ in the center $Z$, take the
localization $A[f^{-1}]$ of the algebra $A$ and glue them together over
the principal open subsets $D(f) = \Spec Z[f^{-1}] \subset \Spec Z$.
\end{proof}

\begin{sloppypar}
\begin{ex}[The quantum space $\nc^n_\epsilon$ at a
root $\epsilon$ of unity]
Con\-sid\-er the associative algebra $A = \nc \langle x_1, \dots , x_n
\rangle / (x_i x_j - \epsilon x_j x_i) \; | \; i < j)$, where $\epsilon
\in \nc$, $\epsilon^l = 1$. The center of $A$ is the subalgebra $Z$
generated by $x_1^l, \dots , x_n^l$, and the whole algebra $A$ is a
free $Z$-module of rank $l^n$.
\end{ex}
\end{sloppypar}

\begin{sloppypar}
\begin{ex}[The quantum group $\SL_{\epsilon} (n)$ at a root of unity]
This \linebreak[0] is a group of linear transformations of
$\nc^n_\epsilon$, see Manin \cite{manin:montreal} or
Faddeev-Reshetikhin-Takhtajan \cite{frt}. To avoid cumbersome
notation, consider the case $n=2$. The quantum group $\SL_q (2)$ is
the associative algebra $A =\nc[q,q^{-1}] \langle a,b,c,d \rangle/ (ab
= q^{-1} ba, \; ac = q^{-1} ca,\; cd = q^{-1} dc, \; bd = q^{-1} db,
\; bc = cb, \; ad -q^{-1} bc = da - q cb = 1)$. Of course, it has a
quantum group, i.e., a Hopf algebra structure coming from the matrix
multiplication of the variables arranged in a matrix $M =
\begin{pmatrix} a & b\\c & d\end{pmatrix}$, which is of no interest so
far, as long as we are concerned only with the noncommutative space
properties of $\SL_q (2)$. When $q= \epsilon$ is an $l$th root of 1,
the subalgebra generated by the entries of $M^l$ is central and the
whole algebra is finite-dimensional over it. The noncommutative space
is then the sheaf corresponding to the algebra $A$ over the usual
algebraic group $\SL(2)$.  This example is in a certain sense dual to
the quantum algebra $U_{\epsilon}({\mathfrak
sl}_2)$ of Example \ref{kac}.
\end{ex}
\end{sloppypar}

\section{Quantum projective spaces}

Let $S = \bigoplus_{n \ge 0} S_n$ be a $\nz$-graded associative
algebra, $Z$ a graded central subalgebra of it, such that $S$ is
finitely generated as a $Z$-module. Then one can construct the {\it
noncommutative projective spectrum} $\Proj_Z S$ as the ringed space
$(X, \Q)$, where $X = \Proj Z$ and the sheaf $\Q$ is constructed as in
Lemma~\ref{lem}: for every homogeneous element $f \in Z$, consider the
localization $S_{(f)} = S[f^{-1}]_0$ and glue them over open subsets
$D((f)) = \Spec Z[f^{-1}]_0$. For a finite-type graded module $M$ over
$S$, we can similarly construct a ``coherent'' sheaf $\mathcal M$ of
$\Q$-modules over $X$ and define the {\it functor of global sections}
\[
\Gamma (X, \mathcal M) = H^0 (X, \mathcal M) = \Hom_{\Q} (\Q, \mathcal M) =
\Hom_{\OO} (\OO, \mathcal M)
\]
just as usual and the {\it cohomology} as its derived functors or as
\[
H^q (X, \mathcal M) = \Ext^q_{\Q} (\Q, \mathcal M) = \Ext^q_{\OO} (\OO,
\mathcal M),
\]
where $\OO$ is the structure sheaf of $X = \Proj Z$. We can define an {\it
invertible sheaf} $\Q(n)$, $n \in \nz$, as the sheaf of $\Q$-modules
associated with the graded free $S$-module $S(n)$ of rank one defined
by $S(n)_m = S_{m+n}$.

\subsection{The quantum projective space $\np^n_\epsilon$}

Consider the graded algebra
\begin{equation}
\label{S}
S = \nc_q \langle X_0, \dots , X_n
\rangle/ (X_i X_j = q X_j X_i \text{ for $i < j$}),
\end{equation}
over $\nc_q = \nc[q,q^{-1}]$. We can think of it as the homogeneous
ring of a quantum projective space $\np^n_q$.  When $q$ is specified
to a root $\epsilon$ of unity, $S$ gets a huge center, the subalgebra
$Z$ generated by $X_0^l, \dots , X_n^l$. This subalgebra lifts to a
subspace closed under the deformed multiplication at an arbitrary
value q. It is still the subalgebra generated by the $l$-th powers of
the generators, with the same kind of multiplication law as for the
original projective space, but with $q$ replaced by $q^{l^2}$, which
is not 1 when $q$ is not equal to an $l^2$'th root of unity. Thus, we
are in the situation of a Poisson fibred algebra with zero curvature,
which we studied before. We will save studying those structures for a
future work; for now, we would like to focus on the scalar
multiplication structure, which is already very rich.

The {\it quantum projective space} $X = \np^n_\epsilon$ is the pair $(X,
\Q_X)$, where $X =\np^n$ and the sheaf $\Q_X$ of noncommutative algebras
obtained from the noncommutative graded $\nc[X_0^l, \dots X_n^l]$-algebra
\eqref{S}. The sheaf $\Q_X$ can be described as follows.  Consider the degree
$l^n$ ramified covering
\begin{align*}
\pi: \np^n & \to \np^n, \\
(X_0: \dots : X_n ) & \mapsto (X_0^l: \dots : X_n^l) .
\end{align*}
The sheaf $\bar{\Q}_X = \pi_* \OO$ on $\np^n$ is
\begin{equation}
\label{dec}
\bar{\Q}_X = \OO \oplus m_1 \OO (-1) \oplus m_2 \OO(-2) \oplus \dots \oplus m_n
\OO(-n),
\end{equation}
where $m_1 = p(l,n+1;l) =  \begin{pmatrix}n+l\\n\end{pmatrix} - n - 1$
is the number of partitions of $l$ into $n+1$ nonnegative parts, each
strictly smaller than $l$, $m_2 = p(2l,n+1;l) =
\begin{pmatrix}n+2l\\2l\end{pmatrix} - (n+1) m_1 -
\begin{pmatrix}n+2\\2\end{pmatrix}$, $m_3 = p(3l, n+1;l) =
\begin{pmatrix}n+3l\\3l\end{pmatrix} - (n+1) m_2 -
\begin{pmatrix}n+2\\2\end{pmatrix} m_1 -
\begin{pmatrix}n+3\\3\end{pmatrix}$, \dots, $m_n = p(nl, n+1;l)$.
$\bar{\Q}_X$ is a sheaf of algebras with the product coming from that on
$\OO$ upstairs. It can be described in terms of the decomposition
\eqref{dec}. For instance, when $n=1$, we obtain $\bar{\Q}_X = \OO
\oplus (l-1) \OO(-1)$, the $l-1$ components $\OO(-1)$ being generated
by $X_0 X_1^{l-1}, X_0^2 X_1^{l-2}, \dots, X_0^{l-1} X_1$ in degree 1,
with the following multiplication table on the generators
\begin{multline*}
X_0^p X_1^{l-p} \cdot X_0^q X_1^{l-q} \\
= X_0^{p+q} X_1^{2l-p-q} =
\begin{cases}
X_1^l (X_0^{p+q} X_1^{l-p-q}) & \text{if $p+q < l$,}\\
X_0^l X_1^l \cdot 1 & \text{if $p+q = l$,}\\
X_0^l (X_0^{p+q-l} X_1^{2l-p-q}) &  \text{if $p+q > l$,}
\end{cases}
\end{multline*}
which induces the following mappings of sheaves:
\[
\OO(-1) \otimes \OO(-1) \quad
\begin{cases}
\begin{CD}
@>X_1^l \cdot>> & \OO(-1),\\
@>X_0^l X_1^l \cdot>> & \OO,\\
@>X_0^l \cdot>> & \OO(-1),
\end{CD}
\end{cases}
\]
respectively.

Thus when $n=1$, we get
\[
X_0^p X_1^{l-p} \cdot X_0^q X_1^{l-q} = \epsilon^{(p-l)q} X_0^{p+q}
X_1^{2l-p-q}.
\]
Note that this multiplication law is still commutative --- what else can we
expect from a one-dimensional space? The sheaves $Q_X$ of quantum functions on
$\np^n_\epsilon$ will be, of course, noncommutative, starting from $n =2$.

\begin{rem}
One can think of the morphism of ringed spaces
\[
F: \np^n_\epsilon = (X, \Q) \to \np^n = (X, \OO)
\]
defined by the natural embedding
\begin{align*}
\nc[Z_0, \dots, Z_n] & \to S,\\
Z_i & \mapsto X_i^l,
\end{align*}
as a {\it quantum Frobenius morphism}.
\end{rem}

\subsection{Cohomology of invertible sheaves on $\np^n_\epsilon$}

The quantum projective spaces $\np^n_\epsilon$ are the simplest examples of
compact quantum homogeneous spaces $G/P$ at roots of unity. Contrary to
the case of generic $q$, where such spaces were studied at the level of
their noncommutative algebras of functions, see
Lakshmibai-Reshetikhin \cite{lr}, Soibelman \cite{soil}
and Chari-Pressley \cite{chari}, our root-of-unity spaces are
honest topological spaces, ringed with rather small sheaves of
noncommutative algebras. This has the advantage that one can use the
usual cohomology theory of sheaves of abelian groups on them.

As usual, the cohomology of invertible sheaves on $\np^n_\epsilon$
comprise a class of highest-weight representations of the quantum
group $\SL_{\epsilon}(n)$. In case of $\SL_{\epsilon} (2)$, when
$\np^1_\epsilon$ is nothing but the complete flag space $G/B$, we will
obtain the Borel-Weil-Bott theorem, which realizes all irreducible
representations in the cohomology.

First, we should describe the ``Picard group'', which is here nothing
more than a pointed set a priori, of $\np^n_\epsilon$.
\begin{prop}
The Picard set $\Pic \np^n_\epsilon$ of isomorphism classes of
invertible left $\Q$-modules on $\np^n_\epsilon$ is isomorphic to
$\nz$, the sheaves $\Q(j)$, $j \in \nz$, corresponding to the graded
$S$-modules $S(j) = \bigoplus_{m \ge -j} S(j)_m = \bigoplus_{m \ge -j}
S_{j+m}$ making up a complete list of representatives.
\end{prop}

%\begin{sloppypar}
\begin{proof}
As usual, the Picard set can be computed as the nonabelian
cohomology $H^1(\np^n, \Q^*)$, which we will compute here for the
standard covering of $\np^n$ with open sets $U_i = \{X^l_i \ne 0\}$,
$i = 0, 1, \dots, n$, in terms of homogeneous coordinates $X_i$ on
$\np^n_\epsilon$. A cocycle defining a class in $H^1(\np^n, \Q^*)$ can
be represented by functions $f_{ij} \in \Gamma (U_{ij}, \Q^*) \; | \;
0 \le i,j \le n$. Each function $f_{ij}$ is an invertible degree 0
element of $\nc[X_0,\dots, X_n; \linebreak[0] X_i^{-l}, X_j^{-l}]$. Therefore,
$f_{ij} = c_{ij} X_i^d X_j^{-d}$. Note that $d$ is independent of the
choice of a pair $ij$. This defines a natural degree mapping
\begin{equation}
\label{deg}
\Pic \np^n_\epsilon \to \nz,
\end{equation}
which assigns the degree $d$ in $X_i$ of the function $f_{ij}= c_{ij}
X_i^d X_j^{-d}$ representing a cocycle in $H^1(\np^n, \Q^*)$.

We claim that \eqref{deg} is an isomorphism of pointed sets. Indeed, the
inverse of \eqref{deg} is given by assigning the class of the class of
$\Q(d)$ to each $d \in \nz$. It is a true inverse, because any class
$\{f_{ij}\} \in H^1(X, \Q^*)$ is of the form $c_{ij} X_i^d X_j^{-d}$,
with $c_{ij}$ forming a 1-cocycle with coefficients in the constant
sheaf $\nc^*$. Since the cohomology of a noetherian topological space
with constant coefficients is always zero in all degrees higher than
zero, the 1-cocycle $\{c_{ij}\}$ is a coboundary and thus,
$\{f_{ij}\}$ is equivalent to the class of $ \{X_i^d X_j^{-d}\}$,
which is the class of $\Q(d)$ in $\Pic \np^n_\epsilon$.
\end{proof}
%\end{sloppypar}

The quantum group $\SL_{\epsilon} (n+1)$ acts naturally on
$\np^n_\epsilon$ and  the sheaves $\Q(j)$. This means the sheaf $\Q$ of
functions and the sheaves of sections of $\Q(j)$ have a (right) $\SL_\epsilon
(n+1)$-comodule structure, coming from the natural comodule structure
on the graded algebra $S$ and grading shifts of it.

\begin{thm}
\label{coh}
The cohomology groups $H^i(X, \Q(j))$ of the invertible sheaves on the
quantum projective space $X = \np^n_\epsilon$ as representations of
$\SL_{\epsilon} (n+1)$ are given by the following formula:
\[
H^i(X, \Q(j)) =
\begin{cases}
S^j (V^*), & \text{for $j \ge 0$, $i =0$,}\\
S^{-n-1-j} (V)  & \text{for $j \le -n-1$, $i =n$,}\\
0 & \text{otherwise,}
\end{cases}
\]
where $V^*$ is the standard $n+1$-dimensional representation of
$\SL_\epsilon (n+1)$ in the space $V = S_1$, a graded component of the
algebra $S = \nc \langle X_0, \dots, X_n \rangle / (X_i X_j = \epsilon
X_j X_i \text{ for } i < j)$. The $j$th symmetric power is understood
as the representation of $\SL_\epsilon (n+1)$ in $S_j$.
%and the exterior power $\Lambda^{n+1} (V^*)$as the representation in
%the highest component $A_{n+1}$ of the Koszul dual algebra $A = \nc
%\langle \xi_0, \dots, \xi_n \rangle / (\xi^2 = 0, \xi_i \xi_j = -
%\epsilon^{-1} \xi_j \xi_i \text{ for } i < j)$, which is the algebra of
%functions of the quantum superspace $\nc_\epsilon^{0|n+1}$. This
%representation is of course trivial for the group $\SL_\epsilon$.
\end{thm}

\begin{proof}
The space $H^0 (X, \Q(j))$ can be naturally identified with the space of
globally defined quantum rational functions of $X_i$ of degree $j$. These are
just quantum polynomials in $X_i$ of degree $j$, that is
the space $S_j$. The computation of $H^n$ is given by the natural
$\SL_\epsilon$-equivariant pairing  $\Q(j) \otimes \Q(-n-1-j) \to \Q(-n-1)$ and
the trace mapping $H^n (X, \Q (-n-1)) \to \nc$, which follows from the
decomposition
\begin{equation*}
\label{dec'}
\Q (j) = \bigoplus_{k \in \nz} p(kl+j, n+1; l) \OO (-k)
\end{equation*}
(only a finite number of terms do not vanish), generalizing
\eqref{dec}. The same decomposition and the computation of cohomology
of $\OO(j)$ implies vanishing of the other cohomology groups.
\end{proof}

\section{The Borel-Weil-Bott theorem}

Here we are going to apply Theorem~\ref{coh} to deduce
Corollary~\ref{BWB}, the quantum version of the Borel-Weil-Bott (BWB)
theorem in the simplest case of the quantum group $\SL_{\epsilon} (2)$
when $\epsilon$ is an $l$th root of unity for $l$ odd. Thus, it will
complete the program indicated in Manin \cite{manin:montreal} and
started in Parshall-Wang \cite{pw}, Andersen \cite{andersen}, and
Lakshmibai-Reshetikhin \cite{lr} entirely in terms of the functor of
induction from the quantum Borel subgroup, without any use of the
geometry of homogeneous spaces.

We need to recall some basic facts of representation theory of $G =
\SL_{\epsilon} (2)$, see \cite{pw}. There is a natural sequence of
quantum group embeddings, i.e., surjective inverse morphisms of Hopf
algebras,
\[
T \subset B \subset G,
\]
where $T$ is the quantum subgroup of diagonal matrices, which is just
the usual torus $\nc^*$, and $B$ is the subgroup of upper triangular
matrices. Given a character $\lambda \in \nz$ of $T$, we can extend it
trivially to a character of $B$ and then form an induced representation
$\ind_B^G \lambda$, a highest weight module. It is not equal to zero if
and only if $\lambda \ge 0$. If $\lambda < l$, it is irreducible.
Another, purely root-of-unity construction of irreducibles starts from
an irreducible module $L (\lambda)$, $\lambda \ge 0$, of the ordinary
group $\SL(2)$ and pulls it back to a representation of $\SL_\epsilon
(2)$ via the {\it quantum Frobenius morphism}
\[
F: \SL_\epsilon (2) \to \SL (2)
\]
sending the standard generators $A, B, C, D$ of the algebra of regular
functions on $\SL(2)$ to $a^l, b^l, c^l, d^l$ in the Hopf algebra
$\SL_\epsilon (2)$, respectively. The pullback $F^* L (\lambda)$ is an
irreducible representation $L(l \lambda)$ of $G= \SL_\epsilon (2)$, sitting
inside $\ind_B^G l \lambda$. Finally, the complete list of irreducibles is
given by tensor products
\[
L(\lambda_0 + l \lambda_1) = L(\lambda_0) \otimes F^* L (\lambda_1) \subset
\ind_B^G (\lambda_0 + l \lambda_1),
\]
for $0 \le \lambda_0 < l, \quad 0 \le \lambda_1$.

In view of these facts, our computation of cohomology in
Theorem~\ref{coh} yields the following result.
\begin{crl}[Quantum BWB Theorem for $G=\SL_{\epsilon} (2)$]
\label{BWB}

\begin{enumerate}
\item As a representation of $G = \SL_{\epsilon} (2)$,
\[
H^i (X, \Q(\lambda)) =
\begin{cases}
\ind_B^G \lambda & \text{if $ \lambda \ge 0 $ and $i=0$,}\\
\ind_B^G (-2 - \lambda) & \text{if $ \lambda \le -2 $ and $i=1$,}\\
0 & \text{otherwise}.
\end{cases}
\]
\item Any irreducible representation of $G$ is
\begin{align*}
L(\lambda) = H^0 (X, \Q(\lambda)) & \qquad \text{for $0 \le \lambda < l$,}\\
L(l\lambda) = H^0 (X, \OO (\lambda)) & \qquad \text{for $0 \le \lambda$}
\end{align*}
or the image $L(\lambda_0 + l \lambda_1)$ of the natural product
\[
H^0 (X, \Q(\lambda_0)) \otimes H^0 (X, \OO (\lambda_1)) \to H^0 (X,
\Q(\lambda_0 + l \lambda_1))
\]
coming from $\Q(\lambda_0) \otimes \OO (\lambda_1) \to \Q(\lambda_0 + l
\lambda_1)$.
\end{enumerate}
\end{crl}

\bibliographystyle{amsplain}
%\bibliography{quantum}
%\end{document}

\providecommand{\bysame}{\leavevmode\hbox to3em{\hrulefill}\thinspace}

\end{document}